\documentclass{amsproc}

\usepackage{latexsym}
\usepackage{amsmath}
\usepackage{amsfonts}
\usepackage{amssymb}
\usepackage{amscd}

\numberwithin{equation}{section}

\newtheorem{proposition}{Proposition}[section]
\newtheorem{lemma}[proposition]{Lemma}
\newtheorem{theorem}[proposition]{Theorem}

\newtheorem{remark}[proposition]{Remark}

\def\itheorem#1#2{\newtheorem{#1}[proposition]{#2}}

\itheorem{Remark}{Remark}

\def\Hom{\mathop{\rm Hom}\nolimits}
\def\Ext{\mathop{\rm Ext}\nolimits}

\def\cf{\mathop{\rm cf}\nolimits}

\def\Q{{\mathbb{Q}}}
\def\Z{{\mathbb{Z}}}

\begin{document}

\title{A characterization of $\Ext(G,\Z)$ assuming $(V=L)$}

\author{Saharon Shelah}

\address{Department of Mathematics,
The Hebrew University of Jerusalem, Israel, and Rutgers
University, New Brunswick, NJ U.S.A.}

\email{Shelah@math.huji.ac.il}
\thanks{Number 873 in Shelah's list of publications. The first author was supported by
project No. I-706-54.6/2001 of the {\em German-Israeli Foundation
for Scientific Research \& Development}.\\ The second author was
supported by a grant from the German Research Foundation DFG}

\author{Lutz Str\"ungmann}

\address{Department of Mathematics,
University of Duisburg-Essen, 45117 Essen, Germany}

\email{lutz.struengmann@uni-essen.de}

\address{{\it Current address}: Department of Mathematics, University of Hawaii, 2565 McCarthy Mall, Honolulu, HI 96822-2273, USA}

\email{lutz@math.hawaii.edu}

\begin{abstract}
In this paper we complete the characterization of $\Ext(G,\Z)$ for any torsion-free
abelian group $G$ assuming G\"odel's axiom of constructibility plus there is no weakly compact cardinal. In particular, we prove in $(V=L)$ that, for a
singular cardinal $\nu$ of uncountable cofinality which is less
than the first weakly compact cardinal and for every sequence of
cardinals $( \nu_p : p \in \Pi )$ satisfying $\nu_p \leq 2^{\nu}$,
there is a torsion-free abelian group $G$ of size $\nu$ such that
$\nu_p$ equals the $p$-rank of $\Ext(G,\Z)$ for every prime $p$ and
$2^{\nu}$ is the torsion-free rank of $\Ext(G,\Z)$.
\end{abstract}

\maketitle \setcounter{section}{0}

\section{Introduction}
Since the first author solved the well-known Whitehead problem in
1977 (see \cite{Sh1}, \cite{Sh2}) the structure of $\Ext(G,\Z)$ for
torsion-free abelian groups $G$ has received much attention. Easy
arguments show that $\Ext(G,\Z)$ is always a divisible group for
every torsion-free group $G$. Hence it is of the form
\[ \Ext(G,\Z)= \bigoplus\limits_{p \in
\Pi}\Z(p^{\infty})^{(\nu_p)} \oplus \Q^{(\nu_0)} \] for some
cardinals $\nu_p,\nu_0$ ($p \in \Pi)$ which are uniquely
determined. Thus, the obvious question that arises is which
sequences $(\nu_0, \nu_p : p \in \Pi)$ can appear as the cardinal
invariants of $\Ext(G,\Z)$ for some (which) torsion-free abelian
group? On the one hand there are a few results about possible
sequences $(\nu_0, \nu_p : p \in \Pi)$ provable in $ZFC$. For
instance, the trivial sequence consisting of zero entries only can
be realized by any free abelian group. On the other hand, assuming
G\"odel's constructible universe $(V=L)$ plus there is no weakly
compact cardinal it has been shown that almost all sequences (with
natural restrictions) can be the cardinal invariants of
$\Ext(G,\Z)$ for some torsion-free abelian group $G$ whenever the
size of the group $G$ is not a singular cardinal of uncountable
cofinality (see section $2$ for details and \cite{EkHu},
\cite{EkSh}, \cite{GS}, \cite{GS2}, \cite{HHS}, \cite{MRS},
\cite{SS1}, and \cite{SS2} for references). However, the question
of which sequences $(\nu_0, \nu_p : p \in \Pi)$ can occur is
independent of $ZFC$. It is the purpose of this paper to deal with
the remaining case, namely torsion-free abelian groups of
cardinality $\nu$ where $\nu$ is singular of cofinality $\cf(\nu)>
\aleph_0$. The idea is to use the construction principle from
\cite{MRS} which holds under $(V=L)$ and to apply the main
theorem from \cite{MRS} in our construction.\\

Our notation is standard and we write maps from the left. If $H$ is a pure subgroup of the abelian group $G$, then we shall write $H \subseteq_* G$. For
further details on abelian groups we refer to \cite{Fu} and for
set-theoretic methods we refer to \cite{EkMe}, \cite{Je} or
\cite{Ku}.

\section{The structure of $\Ext(G,\Z)$}
In this section we recall the basic results on the structure of
$\Ext(G,\Z)$ for torsion-free $G$. It is easy to see that
$\Ext(G,\Z)$ is divisible for torsion-free $G$, hence it is of the
form \[ \Ext(G,\Z)= \bigoplus\limits_{p \in
\Pi}\Z(p^{\infty})^{(\nu_p)} \oplus \Q^{(\nu_0)} \] for some
cardinals $\nu_p,\nu_0$ ($p \in \Pi)$. Since the cardinals $\nu_p$
($p \in \Pi)$ and $\nu_0$ completely determine the structure of
$\Ext(G,\Z)$ we introduce the following terminology. We denote by
$r^e_0(G)$ the {\it torsion-free rank} $\nu_0$ of $\Ext(G,\Z)$
which is the dimension of $\Q \otimes \Ext(G,\Z)$ and by
$r_p^e(G)$ the $p$-$rank$ $\nu_p$ of $\Ext(G,\Z)$ which is the
dimension of $\Ext(G,\Z)[p]$ as a vector space over $\Z/p\Z$ for
any prime number $p \in \Pi$. There are only a few results
provable in $ZFC$ when $G$ is uncountable, but assuming G\"odel's
universe an almost complete characterization is known (if there is no weakly compact cardinal). The aim of this paper is to fill the remaining gap.\\
We first justify our restriction to torsion-free $G$. Let $A$ be
any abelian group and $t(A)$ its torsion subgroup. Then
$\Hom(t(A),\Z)=0$ and hence we obtain the short exact sequence
\[0 \rightarrow \Ext(A/t(A),\Z) \rightarrow \Ext(A,\Z)
\rightarrow \Ext(t(A),\Z) \rightarrow 0 \] which must split since
$\Ext(A/t(A),\Z)$ is divisible. Thus \[ \Ext(A,\Z) \cong
\Ext(A/t(A),\Z) \oplus \Ext(t(A),\Z). \] Since the structure of
$\Ext(t(A),\Z) \cong \prod_{p \in \Pi} \Hom(A,\Z(p^{\infty}))$ is
well-known in $ZFC$ it is reasonable to assume that $A$ is
torsion-free and, of course, non-free. Using Pontryagin's theorem
one proves

\begin{lemma}
\label{countabletf} Suppose $G$ is a countable torsion-free group
which is not free. Then $r_0^e(G)=2^{\aleph_0}$.
\end{lemma}

\proof See \cite[Theorem XII 4.1]{EkMe}. \qed\\

Similarly, we have for the $p$-ranks of $G$ the following result
due to C.U. Jensen.

\begin{lemma}
\label{countablet} If $G$ is a countable torsion-free group, then
for any prime $p$, either $r_p^e(G)$ is finite or $2^{\aleph_0}$.
\end{lemma}

\proof See \cite[Theorem XII 4.7]{EkMe}. \qed\\

This clarifies the structure of $\Ext(G,\Z)$ for countable
torsion-free groups $G$ since the existence of groups as in Lemma
\ref{countabletf} and Lemma \ref{countablet} follows from Lemma
\ref{existence}. We now turn our attention to uncountable groups
and assume G\"odel's axiom of constructibility. The following is
due to Hiller, Huber and Shelah.

\begin{lemma}[V=L]
\label{uncountabletf} Suppose $G$ is a torsion-free non-free group
and let $B$ be a subgroup of $A$ of minimum cardinality $\nu$
such that $A/B$ is free. Then $r_0^e(G)=2^{\nu}$. In particular,
$r_0^e(G)$is uncountable and $r_0^e(G)=2^{|G|}$ if
$G^*=\Hom(G,\Z)=0$.
\end{lemma}

\proof See \cite[Theorem XII 4.4, Corollary XII 4.5]{EkMe}. \qed\\

Note that the above lemma is not true in $ZFC$ since for any
countable divisible group $D$ it is consistent that there exists
an uncountable torsion-free group $G$ with $\Ext(G,\Z) \cong D$,
hence $r_0^e(G)=1$ is possible taking $D=\Q$ (see Shelah \cite{Sh3}).\\
Again we turn to $p$-ranks. There is a useful characterization of
$r_p^e(G)$ using the exact sequence
\[ 0 \rightarrow \Z \overset{p}\rightarrow \Z \rightarrow \Z/pZ \rightarrow
0.
\]
The induced sequence
\[ \Hom(G,\Z) \overset{\varphi^p}{\rightarrow} \Hom(G,\Z/p\Z)
\rightarrow \Ext(G,\Z) \overset{p_*}{\rightarrow} \Ext(G,\Z) \]
shows that the dimension of
\[ \Hom(G,\Z/p\Z)/\Hom(G,\Z)\varphi^p \]
as a vector space over $\Z/p\Z$ is exactly $r_p^e(G)$.

The following result due to Mekler, Roslanowski and Shelah shows
that under the assumption of $(V=L)$ almost all possibilities for
$r_p^e(G)$ can appear if the group is of regular cardinality.

\begin{lemma}[V=L]
\label{uncountableregular} Let $\nu$ be an uncountable regular
cardinal less than the first weakly compact cardinal. Suppose that
$(\nu_p : p\in \Pi)$ is a sequence of cardinals such that for each
$p$, $0 \leq \nu_p \leq 2^{\nu}$. Then there is an almost-free
group $G$ of cardinality $\nu$ such that $r_0^e(G)=2^{\nu}$ and
for all $p$, $r_p^e(G)=\nu_p$.
\end{lemma}

\proof See \cite[Main Theorem 3.9]{MRS}. \qed\\

On the other hand, if the cardinality of $G$ is singular, then the
following holds which was proved by Grossberg and Shelah.

\begin{lemma}
\label{uncountablesingular} If $\nu$ is a singular strong limit
cardinal of cofinality $\omega$, then there is no torsion-free
group $G$ of cardinality $\nu$ such that $r_p^e(G)=\nu$ for any
prime $p$.
\end{lemma}

\proof See \cite[Theorem 1.0]{GS}. \qed\\

Note that Lemma \ref{uncountableregular} shows that the
restriction in Lemma \ref{uncountablesingular} is the only
restriction for singular strong limit cardinals $\nu$ of
cofinality $\omega$. Namely, if $\rho < \nu$ choose a regular
cardinal $\rho \leq \rho' < \nu$ and apply Lemma
\ref{uncountableregular} to obtain a torsion-free group $G'$ with
$r_p^e(G')=\rho$ and $|G'|=\rho'$. Since $\Ext(-,\Z)$ is a
multiplicative functor we can now easily get a torsion-free group
$G$ from $G'$ with $|G|=\nu$ and $r_p^e(G)=\rho$.\\
Also the case
of weakly compact cardinality was dealt with in \cite{SS1} by
Sageev and Shelah.

\begin{lemma}
If $G$ is a torsion-free group of weakly compact cardinality $\nu$
and $r_p^e(G) \geq \nu$ for some prime $p$, then
$r_p^e(G)=2^{\nu}$.
\end{lemma}

\proof See \cite[Main Theorem]{SS1}. \qed\\

The above results show that under the assumption of $(V=L)$ the
structure of $\Ext(G,\Z)$ for torsion-free groups $G$ of
cardinality $\nu$ is clarified for all cardinals $\nu$ except
$\nu$ is singular but not of cofinality $\omega$. This will be the
subject of the next
section.\\
However, in a particular case, namely when $G^*=\Hom(G,\Z)=0$,
then a complete characterization of $\Ext(G,\Z)$ is known in
G\"odel's universe if there is no weakly compact cardinal. The
following is due to Hiller-Huber-Shelah.

\begin{lemma}
\label{dualt} If $G$ is torsion-free such that $\Hom(G,\Z)=0$,
then for all primes $p$, $r_p^e(G)$ is finite or of the form
$2^{\mu_p}$ for some infinite cardinal $\mu_p \leq |G|$. \end{lemma}

\proof See \cite[Lemma XII 5.2]{EkMe}. \qed\\

Together with Lemma \ref{uncountabletf} and the next result due to
Hiller, Huber and Shelah the characterization is complete if
$\Hom(G,\Z)=0$.

\begin{lemma}
\label{existence} For any cardinal $\nu_0$ of the form
$\nu_0=2^{\mu_0}$ for some infinite $\mu_0$ and any sequence of
cardinals $(\nu_p : p \in \Pi)$ less than or equal to $\nu_0$ such
that each $\nu_p$ is either finite or of the form $2^{\mu_p}$ for
some infinite $\mu_p$ there is a torsion-free group $G$ such that
$\Hom(G,\Z)=0$ and $r_0^e(G)=\nu_0$, $r_p^e(G)=\nu_p$ for all
primes $p \in \Pi$.
\end{lemma}

\proof See \cite[Theorem 3(b)]{HHS}. \qed\\

\section{The singular case}

In this section we prove our main theorem which completes the
characterization of $\Ext(G,\Z)$ for torsion-free groups $G$ under
the assumption of G\"odel's axiom of constructibility plus there is no weakly compact cardinal. The idea of
the proof is as follows: For a singular cardinal $\nu$ of
uncountable cofinality we shall construct a torsion-free abelian
group $G$, of size $\nu$, as the union of pure subgroups
$G_{\alpha}$ such that $G$ has prescribed values for $r_0^e(G)$ and
$r_p^e(G)$ $(p \in \Pi)$. Together with the $G_{\alpha}$'s we also
build homomorphisms $f_{\gamma}^p$ for $\gamma < r_p^e(G)$ such that
no non trivial combination $\sum\limits_{l < n}a_lf_{\gamma_l}^p \in
\Hom(G,\Z/p\Z)$ can be factored by the canonical homomorphism
$\varphi^p:\Hom(G,\Z) \rightarrow \Hom(G,\Z/p\Z)$ to a homomorphism in
$\Hom(G,\Z)$. This is a typical application of the diamond principle
which holds under $(V=L)$ for every regular uncountable cardinal. On the other hand we also need that for
every $f \in \Hom(G,\Z/p\Z)$ there is $f_1=\sum\limits_{l <
n}a_lf_{\gamma_l}^p \in \Hom(G,\Z/p\Z)$ and $f_2 \in \Hom(G,\Z)$
such that $f-f_1=\varphi^p(f_2)$. The two demands seem to be hard to go
together but the principle from \cite{MRS} allows us to carry out
the construction.

\begin{theorem}[V=L]
\label{main} Let $\nu$ be an uncountable singular cardinal of
cofinality $\cf(\nu) > \aleph_0$ which is less than the first
weakly compact cardinal. If $(\nu_p : p \in \Pi)$ is a sequence
of cardinals less than or equal to $2^{\nu}$, then there exists a
torsion-free group $G$ such that
\begin{enumerate}
\item $|G|=\nu$;
\item $r_0^e(G)=2^{\nu}$;
\item $r_p^e(G)=\nu_p$ for all $p \in \Pi$.
\end{enumerate}
\end{theorem}

\proof
First we note that by Lemma \ref{existence} there is for every prime $p \in \Pi$ a torsion-free group $G_p$ such that $\Hom(G_p,\Z)=0$, $|G_p|=\nu$, $r_0^e(G_p)=2^{\nu}$ and $r_p^e(G_p)=2^{\nu}=\nu^+$ but $r_q^e(G_p)=0$ for all $q \not= p$. Since $\Ext(-,\Z)$ commutes with direct sums it therefore suffices to assume $\nu_p \leq \nu$ for all $p \in \Pi$. Let $\kappa=\cf(\nu)$ be the cofinality of $\nu$. Chose a continuous increasing sequence $\left< \mu_{\alpha} : \alpha < \kappa \right>$ such that
\begin{itemize}
\item $\lim_{\alpha < \kappa} \mu_{\alpha} = \nu$
\item If $\alpha$ is a successor ordinal, then $\mu_{\alpha}=\lambda_{\alpha}^+$ for some strong limit singular cardinal $\lambda_{\alpha}>\kappa$ such that $\cf(\lambda_{\alpha})=\aleph_0$ and $\lambda_{\alpha+1} > \mu_{\alpha}$
\end{itemize}
Now, let $S\subseteq \{ \alpha < \kappa : \cf(\alpha)=\aleph_0 \}$ be stationary. Inductively we shall construct a torsion-free group $G=\bigcup\limits_{\alpha < \kappa} G_{\alpha}$ such that the following conditions are satisfied.

\begin{enumerate}
\item $|G_{\alpha}|=\mu_{\alpha}$ and $G_{\alpha}$ is $\mu_0$-free; Moreover, there is no $G' \subseteq G_{\alpha}$ of cardinality less than $\mu_{\alpha}$ such that $G_{\alpha}/G'$ is free;
\item $G_{\beta} \subseteq_* G_{\alpha}$ if $\beta < \alpha$;
\item $\Hom(G_{\alpha},\Z/p\Z)=L_p^{\alpha} \oplus K_p^{\alpha}$ such that $f\restriction_{G_{\alpha}} \in K_p^{\alpha}$ for all $f \in K_p^{\beta}$,
$\alpha < {\beta}$;
\item There are bases $B_p^{\alpha}$ of $K_p^{\alpha}$ as vector spaces over $\Z/p\Z$ and functions $T_{\alpha} : B_p^{\alpha} \longrightarrow \Hom(G_{\alpha},\Z)$ such that
\begin{enumerate}
\item $f\restriction_{G_{\alpha}} \in B_p^{\alpha}$ for all $f \in
B_p^{\beta}$, $\alpha < \beta$;
\item $\varphi^pT_{\alpha}=id_{B_p^{\alpha}}$, where $\varphi^p$ is the
canonical map $\varphi^p : \Hom(G_{\alpha},\Z) \rightarrow
\Hom(G_{\alpha},\Z/p\Z)$;
\item If $\alpha < \beta$ and $f \in B_p^{\beta}$, then $T_{\alpha}(f \restriction_{G_{\alpha}})=T_{\beta}(f) \restriction_{G_{\alpha}}$;
\end{enumerate}
\item There are bases $\left< f_{\gamma}^{\alpha,p} : \gamma < \nu_p \cap \mu_{\alpha} \right>$ of $L_p^{\alpha}$ and $M_{\alpha} \subseteq \Hom(G_{\alpha},\Z)$ (for $\alpha \not= 0$) such that
\begin{enumerate}
\item $f_{\gamma}^{\alpha,p} \subseteq f_{\gamma}^{\beta,p}$ for all $\gamma < \nu_p \cap \mu_{\alpha}$, $\alpha < \beta$;
\item $M_{\alpha}=\bigoplus\limits_{\gamma < \nu_p \cap \mu_{\alpha}, \delta \in [\alpha,\kappa)} \Z h_{\gamma,\delta}^{\alpha,p}$ and $\bar{p}h_{\gamma,\delta}^{\alpha,p}=f_{\gamma}^{\alpha,p}$ for all $\delta \in [\alpha,\kappa)$, where $\bar{p}$ is the canonical map $\bar{p}: \Z \rightarrow \Z/p\Z$;
\item for $\alpha < \beta \leq \delta < \kappa$ and $\gamma < \nu_p \cap \mu_{\alpha}$ we have $h_{\gamma, \delta}^{\alpha,p} \subseteq h_{\gamma, \delta}^{\beta,p}$;
\item if $\alpha \in S$ and $g \in \Hom(G_{\alpha +1},\Z)$ and $\bar{p}g \in L_p^{\alpha+1} \backslash \{0\}$, then $g \in M_{\alpha+1}$;
\item if $\nu_p \cap \mu_{\alpha} \leq \gamma < \nu_p \cap \mu_{\alpha+1}$, then $f_{\gamma}^{\alpha+1,p} \restriction_{G_{\alpha}}=0$ and $h_{\gamma,\rho}^{\alpha+1,p} \restriction_{G_{\alpha}}=0$ for all $\rho \in [\alpha+1,\kappa)$.
\end{enumerate}
\end{enumerate}

We first show that it is sufficient to carry out the inductive construction of $G_{\alpha}$ ($\alpha < \kappa)$. Assume that the torsion-free groups $G_{\alpha}$ ($\alpha < \kappa)$ are constructed satisfying conditions (i) to (v). Put $G=\bigcup\limits_{\alpha < \kappa}G_{\alpha}$.
Then $G$ is a torsion-free group of cardinality $\nu$. Moreover, $r_0^e(G)=2^{\nu}$ follows from
Lemma \ref{uncountabletf} and property (i). It remains to prove that $r_p^e(G)=\nu_p$ for all $p \in \Pi$. Let $p \in \Pi$ and $\alpha < \nu_p$. We define
\[ f_{\alpha}^{\kappa,p} = \bigcup\{ f_{\alpha}^{\gamma,p} : \gamma \in [\delta_{\alpha},\kappa) \}\]
where $\delta_{\alpha}=\min\{ \delta < \kappa : \mu_{\delta} \geq
\alpha \}$. By condition (v)(a) the function $f_{\alpha}^{\kappa,p}$
is a well-defined homomorphism $f_{\alpha}^{\kappa,p} \in
\Hom(G,\Z/p\Z)$ for every $\alpha < \nu_p$. We shall show that
\begin{itemize}
\item  $\{ f_{\alpha}^{\kappa,p} : \alpha < \nu_p \}$ are
linearly independent as elements of $\Hom(G,\Z/p\Z)$;
\item $\{ f_{\alpha}^{\kappa,p} : \alpha < \nu_p \}$ are linearly independent in $\Hom(G,\Z/p\Z)$ modulo $\Hom(G,\Z)\varphi^p$,
i.e. no linear combination of them can be factored by $\bar{p}$ to a
homomorphism from $G$ to $\Z$.
;
\item  $\{ f_{\alpha}^{\kappa,p} : \alpha < \nu_p \}$ together with $\Hom(G,\Z)\varphi^p$ generate $\Hom(G,\Z/p\Z)$.
\end{itemize}
Assume first that
\[ \sum\limits_{\alpha \in E}z_{\alpha}f_{\alpha}^{\kappa,p} =0 \]
for some finite subset $E \subseteq \nu_p$ and elements $z_{\alpha} \in \Z/p\Z$. Then there exists $\beta$ such that $\alpha < \nu_p \cap \mu_{\beta}$ and $f_{\alpha}^{\beta,p}=f_{\alpha}^{\kappa,p}\restriction_{G_{\beta}} \not=0$ for all $\alpha \in E$. Hence
\[ \sum\limits_{\alpha \in E}z_{\alpha}f_{\alpha}^{\beta,p}=\left(\sum\limits_{\alpha \in E}z_{\alpha}f_{\alpha}^{\kappa,p}\right) \restriction_{G_{\beta}}=0. \]
But $\left< f_{\alpha}^{\beta,p} : \alpha < \nu_p \cap \mu_{\beta} \right>$ is a basis of $L_p^{\beta}$ and thus $z_{\alpha}=0$ for all $\alpha \in E$. Therefore, the $f_{\alpha}^{\kappa,p}$'s $(\alpha < \nu_p)$ are linearly independent.\\
Now, assume that there exists a finite linear combination $0
\not=\sum\limits_{\alpha \in E}z_{\alpha}f_{\alpha}^{\kappa,p}$
which can be factored by $\bar{p}$ $(0 \not=z_{\alpha} \in \Z/p\Z$
for all $\alpha \in \Z)$. Hence there is $0 \not=g \in \Hom(G,\Z)$
such that
\[\sum\limits_{\alpha \in E}z_{\alpha}f_{\alpha}^{\kappa,p}=\bar{p}g.\]
Since $E$ is finite, there exists $\beta < \kappa$ such that $\alpha
< \nu_p \cap \mu_{\beta}$ for all $\alpha \in E$. Therefore,
\[\bar{p}g \restriction_{G_{\gamma}} = \sum\limits_{\alpha \in E}z_{\alpha}f_{\alpha}^{\kappa,p} \restriction_{G_{\gamma}} = \sum\limits_{\alpha \in E}z_{\alpha}f_{\alpha}^{\gamma,p} \]
for every $\gamma \in [\beta,\kappa)$. By the linear independence
of the $f_{\alpha}^{\kappa,p}$'s we may assume without loss of
generality that $\bar{p}g \restriction_{G_{\gamma}} \not=0$ for
all $\gamma \in [\beta,\kappa)$ since otherwise we sufficiently
enlarge $\beta$ such that
$f_{\alpha}^{\kappa,p}\restriction_{G_{\beta}} \not=0$ for all
$\alpha \in E$. We conclude that $\bar{p}g
\restriction_{G_{\gamma+1}} \in L_p^{\gamma +1} \backslash \{0\}$
for all $\gamma \geq \beta$ and condition (v)(d) implies that $g
\restriction_{G_{\gamma+1}} \in M_{\gamma+1}$ for all $\gamma \in
S$, $\gamma \geq \beta$. Let $\delta=\gamma+1$ for some $\gamma
\in S$, $\gamma \geq \beta$. Then
\[ g \restriction_{G_{\delta}} = \sum\limits_{k < k_{\delta}}b_{k}^{\delta}h_{\alpha_k^{\delta},j_k^{\delta}}^{\delta,p} \]
with $b_k^{\delta} \in \Z \backslash \{0\}$ and $\alpha_k^{\delta} < \nu_p \cap \mu_{\delta}$, $j_k^{\delta} \in [\delta,\kappa)$. Since the $h_{\beta,r}^{\delta,p}$ form a basis of $M_{\delta}$, this representation is unique.\\
By a pigeon hole argument we may assume that $k_{\delta}=k_0$ and
$b_k^{\delta}=b_k$ for arbitrarily large $\delta=\gamma +1$,
$\beta \leq \gamma \in S$. Let $\gamma_0 \in S$ be sufficiently
large such that $k_0=k_{\delta_0}$ with $\delta_0=\gamma_0 + 1$.
Thus
\[ g \restriction_{G_{\delta_0}} = \left( g \restriction_{G_{\delta}} \right) \restriction_{G_{\delta_0}}=
\sum\limits_{k < k_0} b_k
h_{\alpha^{\delta_0}_k,j_k^{\delta_0}}^{\delta_0,p} =
\sum\limits_{k < k_0} b_k
\left(h_{\alpha^{\delta}_k,j_k^{\delta}}^{\delta,p}
\restriction_{G_{\delta_0}} \right)  \] for all $\delta=\gamma+1$,
$\gamma_0 \leq \gamma \in S$. By condition (v)(e) and the two
compatibility conditions (v)(c) and (v)(d) it easily follows that
the following holds for all $\epsilon < \kappa$
\begin{equation}
\tag{v)(e'} \text{if } \nu_p \cap \mu_{\alpha} \leq \gamma < \nu_p
\cap \mu_{\epsilon+1}, \text{ then }\end{equation}
\[
f_{\gamma}^{\epsilon+1,p} \restriction_{G_{\alpha}}=0 \text{ and }
h_{\gamma,\rho}^{\epsilon+1,p} \restriction_{G_{\alpha}}=0 \text{ for
all } \rho \in [\epsilon+1,\kappa).\] Thus, if $\alpha_k^{\delta} \geq
\nu_p \cap \mu_{\delta_0}$ (and also $\alpha_k^{\delta} < \nu_p \cap
\mu_{\delta}$), then (v)(e') implies that $h_{\alpha_k^{\delta},
j_k^{\delta}}^{\delta,p} \restriction_{G_{\delta_0}} =0$. Note that
$\delta=\gamma +1$ for some $\gamma_0 \leq \gamma \in S$. Hence
$\alpha_k^{\delta} < \nu_p \cap \mu_{\delta_0}$ for all $k < k_0$
and therefore \[h_{\alpha_k^{\delta}, j_k^{\delta}}^{\delta,p}
\restriction_{G_{\delta_0}} = h_{\alpha_k^{\delta},
j_k^{\delta}}^{\delta_0,p}\] for all $j_k^{\delta}$. However, by
uniqueness it follows that
\[ \{ \alpha_k^{\delta_0} : k < k_0 \}=\{ \alpha_k^{\delta} : k < k_0
\} \] and also
\[ \{ j_k^{\delta_0} : k < k_0 \}=\{ j_k^{\delta} : k < k_0 \} \]
for all $\delta$ large enough. This contradicts the fact that
$j_k^{\delta} \in [ \delta, \kappa )$ for all $k < k_0$ and
$\delta$. Therefore we obtain $\nu_p \leq r_p^e(G)$ for all $p \in
\Pi$.\\
It remains to prove that $r_p^e(G)=\nu_p$ for all $p \in \Pi$. For
this it suffices to show that $\left< f_{\alpha}^{\kappa,p} : \alpha
< \nu_p \right>$ generate $\Hom(G,\Z/p\Z)$ modulo
$\Hom(G,\Z)\varphi^p$. Hence, let $0 \not= g \in \Hom(G,\Z/p\Z)$. We have to
prove that there is a finite linear combination $\sum\limits_{\alpha
\in E}z_{\alpha}f_{\alpha}^{\kappa,p}$ with $E \subseteq \nu_p$ and
$0 \not=z_{\alpha} \in \Z/p\Z$ such that
\[ g - \sum\limits_{\alpha
\in E}z_{\alpha}f_{\alpha}^{\kappa,p} = \bar{p}h \] for some $h \in
\Hom(G,\Z)$. Let $g_{\alpha}=g \restriction_{G_{\alpha}}$ for all
$\alpha < \kappa$. Hence, by (iii), there exist $k_{\alpha} \in
K_{\alpha}^p$ and $l_{\alpha} \in L_{\alpha}^p$ such that
$g_{\alpha}=k_{\alpha} + l_{\alpha}$ for every $\alpha < \kappa$
since $g_{\alpha} \in \Hom(G_{\alpha},\Z/p\Z)$. Thus
\[ l_{\alpha} = \sum\limits_{\beta \in E_{\alpha}}
z_{\beta}^{\alpha} f_{\beta}^{\alpha,p} \] for some finite subset
$E_{\alpha} \subseteq \nu_p \cap \mu_{\alpha}$ and $0
\not=z_{\beta}^{\alpha} \in \Z/p\Z$. Since this representation is
unique we may assume by a pigeon hole argument that
$z_{\beta}^{\alpha}=z_{\beta}$ and $E_{\alpha}=E$ are independent of
$\alpha < \kappa$. Note that $f_{\beta}^{\alpha,p}
\restriction_{G_{\alpha'}}=f_{\beta}^{\alpha',p}$ if $\alpha' <
\alpha$ and $\beta < \nu_p \cap \mu_{\alpha'}$ and otherwise $0$ by
$(v)(e')$. We conclude that
\[ \bar{h} = g - \sum\limits_{\beta \in E}
z_{\beta}f_{\beta}^{\kappa,p} \] satisfies $\bar{h}
\restriction_{G_{\alpha}}=k_{\alpha}$ for all $\alpha < \kappa$.
Since $B_p^{\alpha}$ forms a basis of $K_p^{\alpha}$ for all $\alpha
< \kappa$ there is a finite subset $F_{\alpha} \subseteq
B_p^{\alpha}$ and $0\not= w_b^{\alpha} \in \Z/p\Z$ for $b \in
F_{\alpha}$ such that
\[ l_{\alpha} = \sum\limits_{b \in F_{\alpha}} w_b^{\alpha}b \]
for all $\alpha < \kappa$. Again, a pigeon hole argument allows us
to assume that $w_b=w_b^{\alpha}$ and $F=F_{\alpha}$ are independent
of $\alpha$ by uniqueness. Note that $b \restriction_{G_{\beta}} \in
B_p^{\beta}$ if $b \in B_p^{\alpha}$ and $\alpha > \beta$. Putting
\[ h= \sum\limits_{b \in F}w_b\bigcup\limits_{\alpha < \kappa}
T_{\alpha}(b) \] it follows that $h \in \Hom(G,\Z)$ is well-defined
by (iv) and hence $\bar{p}h=\bar{h}$. Therefore,
$g-\sum\limits_{\alpha \in E}z_{\alpha}f_{\alpha}^{\kappa,p} =
\bar{h}$ has a lifting to $h \in \Hom(G,\Z)$. This finishes the
proof and it remains to show that we can carry on the induction as
claimed, i.e. we have to construct groups $G_{\alpha}$ ($\alpha < \kappa$)
such that (i) - (v) are satisfied. We shall distinguish four cases.\\
\underline{Case A:} \quad If $\alpha =0$, let
$G_0=\bigoplus\limits_{\mu_0}\Z$. Moreover, put $L_p^0=\{0\}$ and
$K_p^0=\Hom(G_0,\Z/p\Z)$. Since $G_0$ is free, the existence of
$T_0$ is obvious and choosing $M_0=\{0\}$ all conditions (i) - (v)
are satisfied for
$G_0$.\\
\underline{Case B:} \quad If $\alpha$ is a limit ordinal, then we
let $G_{\alpha}=\bigcup\limits_{\beta < \alpha}G_{\beta}$. As before
we define $\delta_{\gamma}=min \{ \delta < \alpha : \mu_{\delta}
\geq \gamma \}$ for $\gamma < \nu_p \cap \mu_{\alpha}$ and let
\[ f_{\gamma}^{\alpha,p}=\bigcup \{ f_{\gamma}^{\beta,p} : \delta_{\gamma} \leq \beta <
\alpha \} \] and similarly for $\delta \in [\alpha, \kappa)$ and
$\gamma < \nu_p \cap \mu_{\alpha}$ we let
\[ h_{\gamma, \delta}^{\alpha,p}=\bigcup\{h_{\gamma,
\delta}^{\beta,p} : \delta_{\alpha} \leq \beta < \alpha \}. \] By
the continuity conditions (v)(a) and (v)(c) this is well-defined.
Hence, also $M_{\alpha} \subseteq \Hom(G_{\alpha},\Z)$ and
$L_p^{\alpha} \subseteq \Hom(G_{\alpha},\Z/p\Z)$ are defined
canonically. Finally, (iii) and the definition of $L_p^{\alpha}$
induce $K_p^{\alpha}$ as
\[ K_p^{\alpha} = \{ f \in \Hom(G_{\alpha},\Z/p\Z) : f
\restriction_{G_{\beta}} \in K_p^{\beta} \textit{ for all } \beta <
\alpha \}. \] The corresponding set $B_p^{\alpha}=\{ f \in
\Hom(G_{\alpha},\Z/p\Z) : f \restriction_{G_{\beta}} \in B_p^{\beta}
\textit{ for all } \beta < \alpha \}$ is a basis for $K_p^{\alpha}$
and the continuity condition (iv)(c) allows to define
$T_{\alpha}=\bigcup\limits_{\beta < \alpha}T_{\beta}$. It is easy to
check that (i)-(v) are now satisfied.\\
\underline{Case C:} \quad If $\alpha = \beta +1$ and $\beta \not\in
S$, then we let $G_{\alpha}=G_{\beta} \oplus
\bigoplus\limits_{\mu_{\alpha}}\Z$. In the obvious way we define
$L_p^{\alpha}$, $K_p^{\alpha}$, $B_p^{\alpha}$, $T_{\alpha}$,
$f_{\gamma}^{\alpha,p}$ for $\gamma < \nu_p \cap \mu_{\alpha}$,
$M_{\alpha}$, and $h_{\gamma, \delta}^{\alpha,p}$ for $\gamma <
\nu_p
\cap \mu_{\alpha}$ and $\delta \in [\alpha,\kappa)$.\\
\underline{Case D:} \quad If $\alpha=\beta +1$ and $\beta \in S$,
then we imitate the proof of the Main Theorem 3.9 of \cite{MRS}.
We would like to avoid repeating the technical and lengthy
construction from \cite{MRS} but instead point out the main
changes for the convenience of the reader. It is then
straightforward to modify the proof of Main Theorem 3.9 and its
main ingredient Theorem 3.4 from \cite{MRS} and to adopt both to
our setting. We are in the following situation: $\lambda_{\alpha}$
is a strong limit singular cardinal strictly greater than
$\kappa=\cf(\nu)$. Moreover,
$\mu_{\alpha}=\lambda_{\alpha}^+=2^{\lambda_{\alpha}}$ is a
regular cardinal, $\cf(\lambda_{\alpha})=\aleph_0$ and
$\lambda_{\alpha}=\lambda_{\beta+1} > \mu_{\beta}$ ($\mu_{\alpha}$
plays the role of $\lambda$ in \cite{MRS}, so it is the successor
of a strong limit singular cardinal). Since we are assuming
$(V=L)$ the prediction principle from \cite{MRS} holds. In Main
Theorem 3.9 from \cite{MRS} it is proved that we can find a
torsion-free group $\tilde{G}$ (denoted by $G$ in \cite[Main
Theorem 3.9]{MRS}) which has prescribed values $\tilde{\nu}_p$ for
$r_p^e(\tilde{G})$. The construction is very similar to our
construction, i.e. homomorphisms $\tilde{f}^p_{\gamma} \in
\Hom(\tilde{G},\Z/p\Z)$ ($\gamma < \tilde{\nu}_p$) are constructed
(denoted by $f_{\lambda}^{p,\xi}$ in \cite[Main Theorem 3.9]{MRS})
which can not be factorized by $\bar{p}$ (see the proof of Main
Theorem 3.9 in \cite[page 346]{MRS}). The main tool is Theorem 3.4
from \cite{MRS} which can be seen as a Step-Lemma since it deals
with the "killing" of only one undesired homomorphism. However, we
are at stage $\alpha$ of our construction, hence we do not want
that our homomorphisms $f_{\gamma}^{\beta,p}$ ($\gamma < \nu_p
\cap \mu_{\beta}$) which we have dealt with so far have no
extension to $G_{\alpha}$, i.e. can not be factorized by $\bar{p}$
but we just require that there are only some extensions, namely a
set of extensions of $f_{\gamma}^{\alpha,p}$ which is assigned to
each $f_{\gamma}^{\alpha,p}$. The role of this set is played by
$\{ h_{\gamma,\delta}^{\alpha,p} : \delta \in [\alpha,\kappa) \}$
(see (v)(b)). Thus the proof of Theorem 3.4 from \cite{MRS}
carries over to our situation (for the case of successor cardinal
of a strong limit singular cardinal of cofinality $\aleph_0$).
Since the role of $\lambda$ in \cite[Theorem 3.4]{MRS} is played
by $\mu_{\alpha}$ in our setting only cases B and C in the proof
of Theorem 3.4 remain. As usual one guesses the undesirable
factorizations and kills them without effecting the work towards
lifting that has been done already. The only difference is that we
do not require that there is no lifting but we allow only the
assigned ones. Now the adjusted version of Theorem 3.4 from
\cite{MRS} is used in the Main Theorem 3.9 from \cite{MRS} in the
case of $\lambda$ being a successor of a strong limit cardinal.
The resulting group serves as our $G_{\alpha}$. \qed

\begin{remark}
We would like to remark that the only reason for the choice of
$\mu_{\alpha}$ as successor of a strong limit singular cardinal of
cofinality $\aleph_0$ (if $\alpha$ is a successor ordinal) is that
this is the easiest situation in the proof of \cite[Main Theorem
3.9]{MRS}. However, the strategy described in Case D of the above
proof of Theorem \ref{main} (i.e. not killing all extensions of a
homomorphism but allowing some of them to survive) works for every
regular uncountable cardinal which is not weakly compact, e.g.
$\aleph_1$. For instance it follows easily for $\aleph_1$ from
\cite[Theorem XII 4.10]{EkMe} using \cite[Lemma XII 4.8 and Lemma
XII 4.9]{EkMe}.
\end{remark}
\goodbreak

\bigskip

\end{document}